\magnification\magstep1
\font\tengoth=eufm10  \font\ninegoth=eufm9
\font\eightgoth=eufm8 \font\sevengoth=eufm7
\font\sixgoth=eufm6   \font\fivegoth=eufm5
\newfam\gothfam \def\goth{\fam\gothfam\tengoth}
\textfont\gothfam=\tengoth
\scriptfont\gothfam=\sevengoth
\scriptscriptfont\gothfam=\fivegoth

\font\ninerm=cmr9  \font\eightrm=cmr8  \font\sixrm=cmr6
\font\ninei=cmmi9  \font\eighti=cmmi8  \font\sixi=cmmi6
\font\ninesy=cmsy9 \font\eightsy=cmsy8 \font\sixsy=cmsy6
\font\ninebf=cmbx9 \font\eightbf=cmbx8 \font\sixbf=cmbx6
\font\nineit=cmti9 \font\eightit=cmti8
\font\ninett=cmtt9 \font\eighttt=cmtt8
\font\ninesl=cmsl9 \font\eightsl=cmsl8
\newskip\ttglue

\def\eightpoint{\def\rm{\fam0\eightrm}
  \textfont0=\eightrm \scriptfont0=\sixrm
  \scriptscriptfont0\fiverm
  \textfont1=\eighti \scriptfont1=\sixi
  \scriptscriptfont1\fivei
  \textfont2=\eightsy \scriptfont2=\sixsy
  \scriptscriptfont2\fivesy
  \textfont3=\tenex \scriptfont3=\tenex
  \scriptscriptfont3\tenex
  \textfont\itfam=\eightit\def\it{\fam\itfam\eightit}%
  \textfont\slfam=\eightsl\def\sl{\fam\slfam\eightsl}%
  \textfont\ttfam=\eighttt\def\tt{\fam\ttfam\eighttt}%
  \textfont\gothfam=\eightgoth\scriptfont\gothfam=\sixgoth
  \scriptscriptfont\gothfam=\fivegoth
  \def\goth{\fam\gothfam\tengoth}
  \textfont\bffam=\eightbf\scriptfont\bffam=\sixbf
  \scriptscriptfont\bffam=\fivebf
  \def\bf{\fam\bffam\eightbf}%
  \tt\ttglue=.5em plus.25em minus.15em
  \normalbaselineskip=9pt \setbox\strutbox\hbox{\vrule
  height7pt depth2pt width0pt}%
  \let\big=\eightbig\normalbaselines\rm}

\def\ninepoint{\def\rm{\fam0\ninerm}
  \textfont0=\ninerm \scriptfont0=\sixrm
  \scriptscriptfont0\fiverm
  \textfont1=\ninei \scriptfont1=\sixi
  \scriptscriptfont1\fivei
  \textfont2=\ninesy \scriptfont2=\sixsy
  \scriptscriptfont2\fivesy
  \textfont3=\tenex \scriptfont3=\tenex
  \scriptscriptfont3\tenex
  \textfont\itfam=\nineit\def\it{\fam\itfam\nineit}%
  \textfont\slfam=\ninesl\def\sl{\fam\slfam\ninesl}%
  \textfont\ttfam=\ninett\def\tt{\fam\ttfam\ninett}%
  \textfont\gothfam=\ninegoth\scriptfont\gothfam=\sixgoth
  \scriptscriptfont\gothfam=\fivegoth
  \def\goth{\fam\gothfam\tengoth}
  \textfont\bffam=\ninebf\scriptfont\bffam=\sixbf
  \scriptscriptfont\bffam=\fivebf
  \def\bf{\fam\bffam\ninebf}%
  \tt\ttglue=.5em plus.25em minus.15em
  \normalbaselineskip=11pt \setbox\strutbox\hbox{\vrule
  height8pt depth3pt width0pt}%
  \let\big=\ninebig\normalbaselines\rm}
\def\bibliography#1\par{\vskip0pt
  plus.3\vsize\penalty-250\vskip0pt
  plus-.3\vsize\bigskip\vskip\parskip
  \message{Bibliography}\leftline{\bf
  Bibliography}\nobreak\smallskip\noindent
  \ninepoint\frenchspacing#1}

\def\ZZ{{\bf Z}}
\def\FF{{\bf F}}
\def\QQ{{\bf Q}}

\font\twelvebf=cmbx12 at 15pt
{\eightpoint
\noindent Preliminary version \hfill Stanford, September, 2009}
\hrule
\bigskip\bigskip
\centerline{\twelvebf Fitting ideals }
\bigskip
\centerline{\twelvebf and the Gorenstein property}
\bigskip\bigskip
\centerline{\bf Burcu Baran}
\smallskip\centerline{\it Department of Mathematics,}
\smallskip\centerline{\it Stanford University,}
\smallskip\centerline{\it Stanford, CA 94305, USA.}
\bigskip\bigskip\bigskip

{\eightpoint
\noindent {\bf Abstract:} Let $p$ be a prime number and $G$ be a finite commutative group such that $p^{2}$ does not divide the order of $G$. In this note we prove that for every finite module $M$ over the group ring $\ZZ_{p}[G]$, the inequality $\#M\,\leq \,\#\ZZ_{p}[G]/{\hbox{Fit}}_{\ZZ_{p}[G]}(M)$ holds. Here, Fit$_{\ZZ_{p}[G]}(M)$ is the $\ZZ_{p}[G]$-Fitting ideal of $M$.}

\beginsection 1. Introduction

Let $R$ be a commutative ring with identity. For a finitely generated $R$-module $M$, we denote the $R$-Fitting ideal of $M$ by Fit$_{R}(M)$. When $R$ is a discrete valuation ring, it is well known that
$${\rm{length}}\,M\,=\,{\rm{length}}\,R/{\rm{Fit}}_{R}(M).$$
for every finite-length $R$-module $M$. (In fact, this is known for one-dimensional local Cohen-Macaulay rings and $M$ of finite length and finite projective dimension; see $[{\bf{2}}$, Lemma $ 21.10.17.3]$ and $[{\bf{4}}$, Thm. $19.1$]). The equality does not hold when $R$ is not a DVR. Indeed, suppose $R$ is a local ring with maximal ideal $\goth m$ and residue field $k$ such that $R$ is not a DVR. Let $M=k \times k$, so the length of $M$ is $2$. The Fitting ideal Fit$_{R}(M)$ is ${\goth m}^{2}$, so the length of $R/{\rm{Fit}}_{R}(M)$ is $1+{\rm{dim}}_{k}\,{\goth m}/{\goth m}^{2}$, which is greater than $2$. Hence, we have a strict inequality
$${\rm{length}}\,M\,< \,{\rm{length}} \, R/{\rm{Fit}}_{R}(M).
$$
Thus, we ask if for certain rings $R$ it is at least true that for every finite-length $R$-module $M$ we have the inequality ${\rm{length}}\,M\,\le \,{\rm{length}} \, R/{\rm{Fit}}_{R}(M)$. Let $p$ be a prime number. We consider this question for $R=A[C]$ where $C$ is a group of prime order $p$ and $A$ is the ring of integers of an unramified finite extension of $\QQ_{p}$. The following is the main result of this paper. It gives an affirmative answer to our question.

\proclaim Theorem 1.1. Let $M$ be a finite $A[C]$-module where $C$ is a group of prime order $p$ and $A$ is   the ring of integers of an unramified finite extension of $\QQ_{p}$. We have the following inequality
$$\#\,M\,\leq \,\# \, A[C]/{\rm{Fit}}_{A[C]}(M).$$
If the ideal Fit$_{A[C]}(M)$ is a principal ideal then we have an equality.

For any finite abelian group $G$ of order not divisible by $p^{2}$, the group ring $\ZZ_{p}[G]$ is a product of rings of the form $A[C]$ as in Theorem $1.1$. If Theorem $1.1$ is true for rings $S$ and $S'$ then it is also true for their direct product $S \times S'$. Hence, the following is a corollary of Theorem $1.1$.

\proclaim Corollary 1.2. Let $p$ be a prime and let $G$ be a finite commutative group for which $p^{2} \not|\, \#\,G$. Then for every finite $\ZZ_{p}[G]$-module $M$ we have the following inequality
$$\#\,M\,\leq \,\# \, \ZZ_{p}[G]/{\rm{Fit}}_{\ZZ_{p}[G]}(M).$$
If the ideal Fit$_{\ZZ_{p}[G]}(M)$ is a principal ideal then we have an equality. 

The local ring $A[C]$ is a Gorenstein ring. This plays a very important role in the proof of Theorem $1.1$. In Section $3$, we prove Proposition $3.6$ relating ideals in $A[C]$ and in its normalization; this is the key proposition. It is an application of a result (Proposition $3.5$) proved in $[{\bf{1}}]$. In the rest of the paper, we use Proposition $3.6$ to exploit the Gorenstein property of $A[C]$. In Section $4$, we fix a short exact sequence $0 \longrightarrow K \longrightarrow A[C]^{t} \longrightarrow M \longrightarrow 0$ for a  finite $A[C]$-module $M$. Let $\FF_{q}$ be the residue field of $A$. We prove that $K/{\goth m}K$ is an $\FF_{q}$-vector space of dimension $t$ if and only if Fit$_{A[C]}(M)$ is a principal ideal. By using that, in Section $5$ we prove the main result (Theorem $1.1$). 
\smallskip
I would like to thank Ren\'e Schoof for his continuous guidance and support.

\beginsection 2. The definition of a Fitting ideal

Let $R$ be a commutative ring with identity and $M$ be a finitely generated $R$-module. Choose a surjective $R$-morphism $f : R^{t} \longrightarrow M $. The $R$-ideal generated by det$(v_{1}$, $v_{2}$, ...,~$v_{t})$, where $v_{1}$, $v_{2}$, ..., $v_{t} \in {\hbox{Ker}}f$, does not depend on $f$ [see ${\bf{3}},\,{\rm{p.}}741]$. It only depends on the $R$-module~$M$.
\bigskip\noindent
{\bf{Definition 2.1.}} The $R$-ideal generated by all det$(v_{1}$, $v_{2}$, ...,~$v_{t})$, where $v_{1}$, $v_{2}$, ..., $v_{t} \in {\hbox{Ker}}f$, is called the $R$-{\it Fitting ideal} of $M$. It is denoted by Fit$_{R}(M)$.
\bigskip\noindent
We have the following proposition.

\proclaim Proposition 2.2. For a finitely generated $R$-module $M$, the following hold.
\smallskip
\item{1.} If $M=R/I$ for an ideal $I$ of $R$, then Fit$_{R}(M)=I$.
\smallskip
\item{2.} If $N$ is another finitely generated $R$-module, then Fit$_{R}(M \times N)={\hbox{Fit}}_{R}(M){\hbox{Fit}}_{R}(N)$.
\smallskip
\item{3.} For any $R$-algebra $B$, we have ${\hbox{Fit}}_{B}(M \otimes_{R} B)={\hbox{Fit}}_{R}(M)B$.

\smallskip\noindent
{\bf{Proof:}} These follow immediately from the definition of a Fitting ideal and properties of the tensor product.
\bigskip\noindent
{\bf{Example 2.3.}} Suppose $L$ is a finitely generated module over a principal ideal domain $D$. Then we have 
$$
L\, \cong\,  \oplus_{i=1}^{t}  D/a_{i}D,
$$
for certain elements $a_{i}$ in $D$. There exists a natural surjective $D$-morphism 
$$
f : D^{t} \longrightarrow L
$$
whose kernel is generated by the vectors $(a_{1},\,0,\,...\,,\,0),\,(0,\,a_{2},\,...\,,\,0),\,...\,,\,(0,\,0,\,...\,,\,a_{t})$.
Therefore, the $D$-ideal Fit$_{D}(L)$ is generated by the product $a_{1}a_{2}\,...\,a_{t}$. With this example we see that if $L$ were to be a finite $D$-module then we would have $\#\,L=\#\,D/{\rm{Fit}}_{D}(L)$.

\beginsection 3. The Gorenstein group ring $A[C]$

In the rest of the paper, we assume that $R=A[C]$ where $C$ is a cyclic group of prime order $p$ and $A$ is the ring of integers of an unramified finite extension of $\QQ_{p}$. Let $\FF_{q}$ be the residue field of $A$, so $q$ is a power of $p$. Suppose $c$ is a generator of $C$. We have the isomorphism
$$
\phi : R \longrightarrow A[T]/((1+T)^{p}-1)
$$
given by $\phi(c)=1+T$. The ring $A[T]/((1+T)^{p}-1)$ is a local ring with maximal ideal $(p,\,T)$ and residue field $\FF_{q}$. As the depth and the Krull dimension of $R$ are both equal to $1$, the local ring $R$ is a Cohen-Macaulay ring. The element $p$ in the maximal ideal of $R$ is an $R$-regular sequence which generates an irreducible ideal in $R$. Therefore, $R$ is a Gorenstein ring. In other words, it has finite injective dimension. In fact, its injective dimension is equal to its Krull dimension which is~$1$.
\medskip\noindent
{\bf{Notation 3.1.}} We denote the unique maximal ideal of $R$ by ${\goth m}$.
\medskip\noindent
{\bf{Remark 3.2.}} The normalization $\widetilde R$ of $R$ in its total quotient ring is~$A \times A[\zeta_{p}]$. Here $\zeta_{p}$ is a primitive $p$-th root of unity. The ring $\widetilde {R}$ is isomorphic to the product $A[T]/(T) \times A[T]/(N)$ where $N={{(1+T)^{p}-1} \over T}$. The ring $\widetilde R$ is a principal ideal ring. We have the short exact sequence
$$0 \longrightarrow R \mathop{\longrightarrow}\limits^{\eta} {\widetilde R} \mathop{\longrightarrow}\limits^{\vartheta} R/{\goth m} \longrightarrow 0,$$
where the map $\eta$ is given by $\eta(r)=(r\,\,{\rm mod}\,T,\,r\,\,{\rm mod}\,N)$ for every $r \in R$, and the map $\vartheta$ is given by $\vartheta(r_{1},r_{2})=r_{1}-r_{2}\,\,{\rm mod}\,{\goth m}$ for every $(r_{1},r_{2}) \in {\widetilde R}$. Thus, the $A$-module ${\widetilde R}/R$ is isomorphic to the residue field $\FF_{q}$ of $R$, and so the quotient $A$-module ${\widetilde R}/R$ has length $1$.  
\medskip\noindent
{\bf{Notation 3.3.}} For any $R$-module $M$, we denote the tensor product $M \otimes_{R} {\widetilde R}$ by $\widetilde M$.
\medskip\noindent
For an $R$-module $M$, there is always the natural $R$-morphism $\psi$ from $M$ to $\widetilde M$ given by    $\psi(m)=m \otimes 1$. We have the following proposition.

\proclaim Proposition 3.4. Let $M$ be an $R$-module and ${\goth m}$ be the maximal ideal of $R$. Consider the natural $R$-morphism $\psi : M \longrightarrow {\widetilde M}$. The cokernel of $\psi$ is isomorphic to $M/{\goth m}M$ through the map $\tau$ given by $\tau(m \otimes (\lambda,\mu))=(\lambda-\mu)m\,\,{\rm mod}\,{\goth m}M$, for every $m \otimes (\lambda,\mu) \in {\widetilde M}$. The kernel of $\psi$ is killed by $\goth m$, so if $M$ is $\ZZ_{p}$-torsion free then the map $\psi$ is injective.

\smallskip\noindent
{\bf{Proof:}} While proving this proposition, to make the computations easy, we identify the ring $R$ with $A[T]/((1+T)^{p}-1)$ via the isomorphism $\phi$ above. Thus, the maximal ideal $\goth m$ of $R$ is $(p, T)$ and ${\widetilde R}$ is equal to $A[T]/(T) \times A[T]/(N)$ where $N={{(1+T)^{p}-1} \over T}$. Consider the short exact sequence in Remark $3.2$. Tensoring this short exact sequence over $A$ with the $R$-module $M$, we obtain the exact sequence
$$0 \longrightarrow K \longrightarrow M  \mathop{\longrightarrow}\limits^{\psi}  {\widetilde M}  \mathop{\longrightarrow}\limits^{\tau} M/{\goth m}M \longrightarrow 0.$$
Since we identified $\tilde R$ with $A[T]/(T) \times A[T]/(N)$, we also identify $\widetilde M$ with $M/TM \times M/NM$. In this exact sequence, the map $\psi$ is given by $\psi(m)=(m\,\,{\rm mod}\,T,\,m\,\,{\rm mod}\,N)$ for every $m \in M$, and the map $\tau$ is given by $\tau(m_{1},m_{2})=m_{1}-m_{2}\,\,{\rm mod}\,{\goth m}M$ for every $(m_{1},m_{2}) \in {\widetilde M}$. The map $\tau$ that we defined here coincides with the map $\tau$ that we defined in the proposition by the identification of $\widetilde M$ with $M/TM \times M/NM$. With this exact sequence it is clear that the cokernel of $\psi$ is isomorphic to $M/{\goth m}M$ through the map $\tau$. Now consider the kernel $K$ of $\psi$ in the above exact sequence. The $R$-module $K$ is equal to $TM \cap NM$, so it is killed by the ideal $(N,T)$. Since $p \in (N,T)$ and the ideal $(p, T)$ is the maximal ideal, we have $(N,T)=(p,T)$. It follows that $K$ is killed by the maximal ideal $\goth m$ of $R$, and in particular by $p$. Thus, if $M$ is a $\ZZ_{p}$-torsion free $R$-module then $K=0$. Hence, we proved the proposition.
\bigskip\noindent
Consider the following proposition concerning general Gorenstein orders over principal ideal domains.

\proclaim Proposition 3.5. Let $\goth O$ be an order over a principal ideal domain. Then the following properties are equivalent:
\smallskip
\item{$-$} $\goth O$ is Gorenstein,
\item{$-$} for any fractional $\goth O$-ideal $\goth a$, we have $({\goth a}:{\goth a}):=\{r \in {\widetilde {\goth O}}:\,\, r{\goth a} \subset {\goth a} \}$ is equal to ${\goth O}$ if and only if $\goth a$ is invertible. Here, $\widetilde {\goth O}$ is the normalization of $\goth O$ in its total quotient ring.

\smallskip\noindent
{\bf{Proof:}} This is Proposition $2.7$ in $[{\bf{1}}]$.
\bigskip\noindent
Let $J$ be any ideal of $R$, then $(J:J):=\{r \in {\widetilde R}:\,\, rJ \subset J \}$ is a ring and we have
$$
R \subset (J:J) \subset {\widetilde R}.
$$
Since the quotient $A$-module ${\widetilde R}/R$ has length $1$ (see Remark $3.2$), the ring $(J:J)$ is equal to either $R$ or ${\widetilde R}$. In the sequel, to prove the main theorem in Section $5$, we will use the following proposition very often to exploit the fact that $R$ is a Gorenstein ring.

\proclaim Proposition 3.6. If the ideal $J$ of $R$ is not a principal ideal, then it is also an ${\widetilde R}$-ideal.

\smallskip\noindent
{\bf{Proof:}} Suppose $J$ is an ideal of $R$ which is not a principal ideal. By the above explanation, $(J:J)$ is either $R$ or $\widetilde R$. Suppose it is equal to $R$. Since $R$ is an order over the principal ideal domain $A$, we use Proposition $3.5$ and we obtain that the $R$-ideal $J$ is invertible. Since $R$ is a local ring, this occurs only when $J$ is a principal ideal generated by an element which is not a zero-divisor. But this contradicts our assumption that $J$ is not a principal ideal. Thus, we have $(J:J)={\widetilde R}$. Hence, $J$ is also an $\widetilde R$-ideal. 
\bigskip\noindent

\beginsection 4. Finite modules over $A[C]$

In this section, we prove some propositions which we will use in Section $5$ in our proof of the main theorem. From now on, we assume that $M$ is a finite $R$-module. Recall that $C$ has prime order $p$ and $A$ is the ring of integers of an unramified finite extension of $\QQ_{p}$ with residue field $\FF_{q}$ and that~$R=A[C]$ is a local ring with unique maximal ideal $\goth m$ and residue field $\FF_{q}$. We fix a short exact sequence
$$
0 \longrightarrow K \longrightarrow R^{t} \longrightarrow M \longrightarrow 0 \eqno(4.1)
$$
of $R$-modules.

\proclaim Proposition 4.1. Consider the $R$-module $K$ in the short exact sequence $(4.1)$. The quotient $K/{\goth m}K$ is an $\FF_{q}$-vector space. We have
 $${\hbox {dim}}_{\FF_{q}}(K/{\goth m}K) \ge t,$$
 with equality holding if and only if $K$ is $R$-free of rank~$t$.

\smallskip\noindent
{\bf{Proof:}} Since the residue field of the local ring $R$ is $\FF_{q}$, the quotient $K/{\goth m}K$ is an $\FF_{q}$-vector space. Let $d={\hbox {dim}}_{\FF_{q}}(K/{\goth m}K)$, so $K$ admits $d$ generators as an $R$-module, by Nakayama's Lemma. By choosing a surjective map $\varphi : R^{d} \longrightarrow K$, we get an exact sequence 
$$ R^{d}  \mathop{\longrightarrow}\limits^{\varphi'}  R^{t} \longrightarrow M \longrightarrow 0.$$
We tensor this exact sequence over $A$ with $A[1/p]=F$. Since $M$ is a finite $R$-module and $R$ is a free $A$-module of rank $p$, we obtain a surjection
$$ (F)^{pd}  \mathop{\longrightarrow}\limits^{\tilde\varphi'} (F)^{pt} \longrightarrow 0.
$$
Hence, this shows that $d \ge t$. Now, suppose $d=t$. Then the surjection $\tilde\varphi'$ is an isomorphism, implying that ${\rm Ker}\,\varphi' \otimes_{A} F=0$. Since ${\rm Ker}\,\varphi' \subset R^{d}$, it does not have nonzero $A$-torsion. This shows that ${\rm Ker}\,\varphi'=0$, and so ${\rm Ker}\,\varphi=0$. Therefore, the map $\varphi$ is an isomorphism, implying that $K$ is $R$-free of rank $t$. Hence, the proposition follows.
\bigskip\noindent

\proclaim Proposition 4.2. Let $M$ be a finite $R$-module. Consider the short exact sequence $(4.1)$. The $R$-module $K$ is free if and only if Fit$_{R}(M)$ is a principal ideal of $R$.

\smallskip\noindent
{\bf{Proof:}} Suppose $K$ is a free $R$-module. Since $M$ is a finite $R$-module, the rank of $K$ is equal to $t$. Then, by definition of the Fitting ideal, the $R$-ideal Fit$_{R}(M)$ is generated by the determinant of the map from $K$ to $R^{t}$, implying that Fit$_{R}(M)$ is a principal ideal. Now, assume that Fit$_{R}(M)$ is a principal ideal of~$R$. Let Fit$_{R}(M)=\alpha R$ where $\alpha \in R$. Note that $\alpha \in R[1/p]^{\times}$ since $M$ is finite. Thus, $\alpha$ is not a zero-divisor in $R$. We claim that there exist $v_{1}$, $v_{2}$, ..., $v_{t} \in K$ such that ~${\rm{det}}(v_{1},\,v_{2},\,...,\,v_{t})=\alpha u$, where $u$ is a unit in $R$. If this were not to be the case, then for every $w_{1}$, $w_{2}$, ..., $w_{t} \in K$  we would have
$$
 ({\hbox{det}}(w_{1},\,w_{2},\,...,\,w_{t})) \subset \alpha {\goth m}
 $$
where ${\goth m}$ is the unique maximal ideal of $R$. Then we would have Fit$_{R}(M) \subset {\goth m}{\hbox{Fit}}_{R}(M)$, so Nakayama's Lemma would imply that Fit$_{R}(M)=0$. But this would contradict with the fact that $M$ is finite. 
\smallskip
Let $r$ be any element of $K$. We solve the linear system
$$
\lambda_{1}v_{1}+\lambda_{2}v_{2}+\,...\,+\lambda_{t}v_{t}=r
$$
with Cramer's Rule. We get that $\lambda_{i}=\delta_{i}/\alpha$ for some $\delta_{i} \in {\hbox{Fit}}_{R}(M)$. Thus, in particular, all $\lambda_{i}$'s are in $R$. This shows that the vectors $v_{1}$, $v_{2}$, ..., $v_{t}$ generate $K$ over $R$. Now suppose
$$
\lambda_{1}v_{1}+\lambda_{2}v_{2}+\,...\,+\lambda_{t}v_{t}=0,
$$
for $\lambda_{i}$'s which are not all zero. This implies that det($v_{1}$, $v_{2}$, ..., $v_{t})=\alpha u=0$, so that $\alpha=0$. This again contradicts with the finiteness of $M$. Thus, all $\lambda_{i}$'s are zero. As a result, we proved that $K$ is a free $R$-module of rank $t$. Hence, the proposition follows.
\bigskip\noindent
We tensor the short exact sequence $(4.1)$ over $R$ with $\widetilde R$ and we obtain the following commutative diagram:
$$\def\normalbaselines{\baselineskip20pt\lineskip3pt
\lineskiplimit3pt}
\matrix{
0& \longrightarrow &K&\longrightarrow & R^{t} & \longrightarrow & M & \longrightarrow &0 \cr
&&\big\downarrow &&\big\downarrow&&\big\downarrow \cr
&&{\widetilde K}&\longrightarrow& {\widetilde R}^{t} & \longrightarrow& {\widetilde M}& \longrightarrow &0 \cr}
$$
Let $H$ be the image of $\widetilde K$ in ${\widetilde R}^{t}$. We have $H={\widetilde R}K$ inside ${\widetilde R}^{t}$, and we also have the commutative diagram of exact sequences. 
$$\def\normalbaselines{\baselineskip20pt\lineskip3pt
\lineskiplimit3pt}
\matrix{
0& \longrightarrow &K&\longrightarrow & R^{t} & \longrightarrow & M & \longrightarrow &0 \cr
&&\big\downarrow &&\big\downarrow&&\big\downarrow \cr
0&\longrightarrow&H&\longrightarrow& {\widetilde R}^{t} & \longrightarrow& {\widetilde M}& \longrightarrow &0\,. \cr} \eqno (4.2)
$$
\smallskip

\proclaim Proposition 4.3. Consider the $R$-modules $K$ and $H$ in the commutative diagram $(4.2)$. The $R$-module $H/K$ is killed by the maximal ideal $\goth m$ of $R$.

\smallskip\noindent
{\bf{Proof:}} Since $H={\widetilde R}K$, we have ${\goth m}H={\goth m}{\widetilde R}K$. As ${\widetilde R}/R$ is isomorphic to $R/{\goth m}$, it follows that ${\goth m}{\widetilde R} \subset R$, and so ${\goth m}H \subset K$. Hence, ${\goth m}$ kills the quotient $H/K$.

\proclaim Proposition 4.4. Consider the $\widetilde R$-module $H$ in the commutative diagram $(4.2)$. It is free of rank $t$, and $H/{\goth m}H$ is an $\FF_{q}$-vector space of dimension $2t$.

\smallskip\noindent
{\bf{Proof:}} In the commutative diagram $(4.2)$, we see that $H$ is a $\widetilde R$-submodule of ${\widetilde R}^{t}$. Since $\widetilde R$ is a product of two discrete valuation rings and the quotient ${\widetilde R}^{t}/H$ is isomorphic to the finite $\widetilde R$-module $\widetilde M$, the $\widetilde R$-module $H$ is free of rank~$t$. Hence, $H$ is isomorphic to ${\widetilde R}^{t}$. Since the residue field of $R$ is $\FF_{q}$ and the quotient ${\widetilde R}/R$ is $\FF_{q}$, we have ${\widetilde R}/{\goth m} = \FF_{q} \times \FF_{q}$. Here we use that the maximal ideal ${\goth m}$ of $R$ is also an $\widetilde R$-ideal (by Proposition $3.6$). Therefore, the quotient $H/{\goth m}H$ is isomorphic to ${\widetilde R}^{t}/{\goth m}{\widetilde R}^{t}$ which is an $\FF_{q}$-vector space of dimension $2t$. Hence, we proved the proposition.

\beginsection 5. The main theorem

Recall that $M$ is a finite $R$-module and we have the short exact sequence $(4.1)$. In this section, our aim is to prove the following theorem. 

\proclaim Theorem 5.1. Let $M$ be a finite $R$-module. We have 
$$
\# M\,\,\le\,\,\# R/{\hbox{Fit}}_{R}(M),
$$ 
with equality holding when ${\hbox {dim}}_{\FF_{q}}(K/{\goth m}K) = t$.

\medskip\noindent
{\bf{Proof:}} Consider the short exact sequence $(4.1)$. In Proposition $4.1$ we proved that ${\hbox {dim}}_{\FF_{q}}(K/{\goth m}K) \ge t$. Thus, we split the proof of this theorem in two cases.
\smallskip\noindent
{\it Case 1:} Suppose ${\hbox {dim}}_{\FF_{q}}(K/{\goth m}K)=t$.
\smallskip\noindent
Consider the short exact sequence $(4.1)$ for a finite $R$-module $M$. By Proposition $4.1$ the $R$-module $K$ is free of rank $t$. Tensoring the short exact sequence $(4.1)$ over $R$ with $\widetilde R$, we obtain the following commutative diagram:
$$
\def\normalbaselines{\baselineskip20pt\lineskip3pt
\lineskiplimit3pt}
\matrix{0&\longrightarrow&0& \longrightarrow &R^{t}&\mathop{\longrightarrow}\limits^{\varphi} & R^{t} & \longrightarrow & M & \longrightarrow &0 \cr
&&\big\downarrow&&\big\downarrow &&\big\downarrow&& \big\downarrow\rlap{$\scriptstyle \psi$}\cr
0&\longrightarrow &{\rm Ker}\,\widetilde\varphi& \longrightarrow &{\widetilde R}^{t}&\mathop{\longrightarrow}\limits^{\widetilde \varphi}& {\widetilde R}^{t} & \longrightarrow& {\widetilde M}& \longrightarrow &0 \cr
}
$$
Consider the bottom exact sequence. We tensor it over $A$ with $F:=A[1/p]$. Since $F$ is $A$-flat, we obtain an isomorphism
$$0 \longrightarrow {F}^{pt} \mathop{\longrightarrow}\limits^{\widetilde{\widetilde \varphi}} {F}^{pt} \longrightarrow 0.$$
As ${\rm{Ker}}\,{\widetilde\varphi} \subset {\widetilde R}^{t}$, it does not have $A$-torsion. Hence, we have ${\rm{Ker}}\,{\widetilde\varphi}=0$. Now, we apply the snake lemma to this commutative diagram. Since ${\widetilde R}^{t}/R^{t}$ is $\FF_{q}^{t}$, we see that
$
\#\, {\rm Ker}\,\psi= \#\,{\rm Coker}\,\psi$. This implies that
$$
\#\, M= \#\,{\widetilde M}.
$$

Now consider the short exact sequence
$$0\longrightarrow R \longrightarrow R \longrightarrow R/ ({\rm det}\,\varphi)\,R \longrightarrow 0$$
where the first map is given by multiplying by ${\rm det}\,\varphi$ and the second map is the natural quotient map. Since ${\rm det}\,\varphi$ is not a zero-divisor, the quotient $R/ ({\rm det}\,\varphi)\,R$ is finite. In the same way, we tensor this short exact sequence over $R$ with $\tilde R$. Then, we obtain a commutative diagram to which we also apply the snake lemma and get
$$
\#\, R/({\hbox{det}}\,\varphi) \,R= \# \,{\widetilde R}/({\hbox{det}}\,\varphi) \,{\widetilde R}.
$$
The $R$-Fitting ideal of $M$ is generated by ${\rm det}\,\varphi$. By Proposition $2.2(3)$, the $\widetilde R$-Fitting ideal of $\widetilde M$ is also generated by ${\rm det}\,\varphi$. Since $\widetilde R$ is the product of principal ideal domains, we have the equality $\#{\widetilde M}= \#{\widetilde R}/$Fit$_{\widetilde R}(\widetilde M)$ (see Example $2.3$). Therefore, the equality $\#M=  \#R/$Fit$_{R}(M)$ follows, as required.
\smallskip\noindent
{\it Case 2:} Suppose ${\hbox {dim}}_{\FF_{q}}(K/{\goth m}K) > t$.
\smallskip\noindent
 Consider the short exact sequence~$(4.1)$. By Proposition $4.1$ the $R$-module $K$ is not free. Hence, by Proposition $4.2$ the $R$-ideal Fit$_{R}(M)$ is not a principal ideal. We have the following equalities.
$$\eqalign{ \# {\widetilde M}\,&=\, \# {\widetilde R}/ {\hbox{Fit}}_{\widetilde R}({\widetilde M})\,\,\,\,\,\,\,\,\,\,\,\,\,\,\,\,\,\,\,\,\,\,\,\,\,\,{\hbox{by Example $2.3$,}}\cr
&=\, \# {\widetilde R}/ {\hbox{Fit}}_{ R}( M){\widetilde R}\,\,\,\,\,\,\,\,\,\,\,\,\,\,\,\,\,\,\,\,\,\,{\hbox{by Proposition $2.2(3)$,}}\cr
&=\, \# {\widetilde R}/ {\hbox{Fit}}_{ R}( M)\,\,\,\,\,\,\,\,\,\,\,\,\,\,\,\,\,\,\,\,\,\,\,\,\,\,\,{\hbox{by Proposition $3.6$,}}\cr
&=\, \# {\widetilde R}/R\cdot \# R/ {\hbox{Fit}}_{ R}( M)\cr
&=q \cdot \# R/ {\hbox{Fit}}_{ R}( M).\cr}$$
Thus, to show that $\#M\le \#R/{\hbox{Fit}}_{R}(M)$, it is enough to show $q\cdot\#M \le  \#{\widetilde M}$. Let $N$ and $N'$ be the finite $R$-modules fitting into an sequence
$$
0 \longrightarrow N \longrightarrow M  \mathop{\longrightarrow}\limits^{\psi}  {\widetilde M}  \longrightarrow  N'\longrightarrow0
$$
with the natural map $\psi$. Applying the snake lemma to $(4.2)$ then yields the exact sequence of $\FF_{q}$-vector spaces (see Proposition $4.3$).
$$
0 \longrightarrow N\longrightarrow H/K  \longrightarrow \FF_{q}^{t}  \longrightarrow  N' \longrightarrow 0.
$$
It follows that
$$\#  {\widetilde M}/ \# M\cdot\#H/K\,=\,q^{t}.$$
Hence, to show that $q \cdot \# M \le \#{\widetilde M}$, we only need to show that $\# H/K < q^{t}$. By Proposition~$4.3$, we have 
$${\goth m}H \subset K \subset H.$$
Since ${\goth m}H={\goth m} {\widetilde R} K$, we have ${\goth m}H={\goth m}K$ by Proposition $3.6$. Thus, we have $[K : {\goth m}H]=[K : {\goth m}K]$, which is greater than $q^{t}$ by assumption. By Proposition $4.4$, the order of $H/{\goth m}H$ is equal to $q^{2t}$. Therefore, the equality
$$
\#H/{\goth m}H\,=\,\#H/K \cdot \#K/{\goth m}K$$
implies that $\#H/K<q^{t}$. Hence, the theorem follows.

\bigskip\bigskip\noindent
$$\bf\rm{REFERENCES}$$

\smallskip
\item{[{\bf{1}}]} Buchmann, A. J. and Lenstra, H. W., {\sl Approximating rings of integers in number fields}, Journal de Th\'eorie des Nombres de Bordeaux {\bf 6} (1994), 221--260.
\smallskip
\item{[{\bf{2}}]} Grothendieck, A. and Dieudonn\'e J., {\sl EGA $IV$4 Etude locale des sch\'emas et des morphismes de sch\'emas}, Publ. Math. IHES {\bf 32} (1967).
\smallskip
\item{[{\bf{3}}]} Lang, S., {\sl Algebra}, Graduate Texts in Mathematics {\bf 211} (2002), Springer.
\smallskip
\item{[{\bf{4}}]} Matsumura, H., {\sl Commutative ring theory}, Cambridge studies in advanced mathematics {\bf 8} (1986), Cambridge University Press.

\bye